\documentclass[12pt]{article}
\usepackage[colorlinks]{hyperref}
\usepackage{color}
\usepackage{graphicx}
\usepackage{graphics}
\usepackage{makeidx}
\usepackage{showidx}
\usepackage{latexsym}
\usepackage{amssymb}
\usepackage{verbatim}
\usepackage{amsmath}
\usepackage{amsthm}
\usepackage{amsfonts}
\usepackage{geometry}

\title{Why Escape Is Faster Than Expected}
\author{H. Attarchi and L.A. Bunimovich}
\begin{document}
	\maketitle
	\noindent
	\begin{abstract}
		We consider chaotic (hyperbolic) dynamical systems which have a generating Markov partition. Then, open dynamical systems are built by making one element of a Markov partition a ``hole" through which orbits escape. We compare various estimates of the escape rate which correspond to a physical picture of leaking in the entire phase space. Moreover, we uncover a reason why the escape rate is faster than expected, which is the convexity of the function defining escape rate. Exact computations are present for the skewed tent map and Arnold's cat map.
		
		\vspace{.5cm}
	\end{abstract}
	%%%%%%%%%%%%%%%%%%%%%%%%%%%%%%%%%%%%%%%%%%%%%%%%%%%%%%%%%%%%%%%%%%%%%%%%%%%%%%%%%%%%%%%%%%%%%%%%%%%%%%%%%%%%%%%
	\section{Introduction}
	Mathematical studies of open dynamical systems started in 1979 \cite{PY}. In this paper, a natural question about what is going to happen if in a billiard table will appear a hole where through which a billiard ball can fall is studied.\par
	Two exact mathematical questions were formulated in this paper. First of all, one should consider a mathematical billiard where a point (mathematical) particle moves. Billiards are Hamiltonian systems and therefore a natural (physical) measure (phase space volume) is preserved under dynamics. But, what if there is a ``hole" of positive measure (otherwise probability to reach a hole is zero)? Would be there then a natural invariant measure for this open system? Such measures are called conditionally invariant measures. They are characterized by the property that at each iterate of the map (at each moment of discrete-time) the same portion of the remaining in the phase space points escapes through the hole. The real interest, though, is only absolutely continuous conditionally invariant measures (a.c.c.i.m.) \cite{DY06}. Another question was about the existence of a natural quantity that characterizes the dynamics of open systems. Such characteristic is the escape rate of orbits through a hole \cite{DY06}.\par
	Of course, physicists studied open systems long before mathematicians \cite{APT13}. It was noticed in real and numerical experiments that escape from chaotic systems is an exponential function of time. Therefore, the factor in front of time in this exponent was naturally called the escape rate. And mathematically escape rate is defined in the same way.\par
	For a long time, the mathematical theory of open systems was dealing essentially with two tasks, which are proving the existence of conditionally invariant measures and the existence of escape rates. A comprehensive description of these efforts is presented in \cite{DY06}. The situation changed when a new natural question was put forth in \cite{BY11}, which asks how the process of escape depends on the position of a hole in phase space? Thanks to emerging research from this question, several fundamental and seemingly obvious previously existing beliefs were (rigorously) proved to be wrong.\par
	First of all, even in the most uniformly hyperbolic (homogeneously expanding distances in the phase space) dynamical systems, escape rate demonstrates strong oscillations when holes are placed at the different parts of the phase space. A natural example is the doubling map of the unit interval $f(x)=2x\ (mod\ 1)$ where at each point (besides $x=1/2$) dynamics expands distances twice. The most striking result in this direction though is that in a typical ergodic system (observe that just ergodic rather than strongly chaotic systems are considered here), there is a continuum of huge ``holes" through which escape rate is arbitrarily small. The sizes (probabilities, measures) of these holes are arbitrarily close to the size of the entire phase space \cite{BV,BY11}. Therefore, the notion of escape rate (although seemingly being so simple and clear) should be taken with great care in theoretical and experimental studies.\par
	Another indication of this was an unexpected observation made in \cite{GDA12} for a doubling map. In this paper, it was demonstrated that the average over elements of a Markov partition of escape rates is larger than expected. It is a surprising result in view of \cite{BY11} where it was shown that for the doubling map escape rate becomes smaller near periodic points. Namely, the escape rate behaves like a constant $C$ times the size of ``hole" when a hole shrinks to a non-periodic (periodic) point, where $C=1$ ($C<1$). This result was generalized to large classes of maps in \cite{FP,KL}.\par
	Having in mind these results, it was natural to assume that a commonly used in physics studies (so-called ``naive") estimate of the escape rate \cite{APT13} should be greater than the average of escape rates taken over elements (separately used as ``holes") of a partition of the phase space. Indeed, a ``naive" estimate is based on the physical picture that the entire phase space is leaking and thus should average out a slow down of escapes at the periodic points which form a negligible (measure zero) subset of the phase space. However, it was shown in \cite{GDA12} numerically that for the doubling map the opposite inequality holds. Thus, contrary to the seemingly natural expectation based on the rigorous results of \cite{BY11}, it turned out that the naive estimate approaches the average of escape rates from below rather than from above in the limit when the size of the hole tends to zero. It is worthwhile to mention that in this limit both these estimates approach each other.\par
	All elements of a natural Markov partition for the doubling map have the same measures. However, if it is not the case, then other candidates for naive estimation of escape rate appear \cite{GDA12}. There was performed a numerical comparison of different naive estimates for the skewed doubling map.\par
	In the present paper, we introduce an estimate for the escape rate, which is applicable for any chaotic dynamical system with generating Markov partition. We prove that the average of escape rates over elements partition indeed always exceed a naive estimation of the escape rate. Most importantly, we uncover a reason for the validity of this inequality, i.e. why ``escape rate is faster than expected". This reason is the convexity of the logarithm function involved in the definition of the escape rate.\par
	Besides, exact numerical examples are presented for the skewed tent map and Arnold's cat map. We also prove for general maps inequality between two naive estimates of the escape rate considered in \cite{GDA12} for the skewed doubling map.
	%%%%%%%%%%%%%%%%%%%%%%%%%%%%%%%%%%%%%%%%%%%%%%%%%%%%%%%%%%%%%%%%%%%%%%%%%%%%%%%%%%%%%%%%%%%%%%%%%%%%%%%%%%%%%%%
	\section{Open dynamical systems and escape rate}
	Assume a map $\hat{T}:\hat{M}\longrightarrow\hat{M}$ generates a dynamical system on the phase space $\hat{M}$ and $H$ be a subset of $\hat{M}$. Then, the open dynamical system corresponding to the map $\hat{T}$ with the hole $H$ is defined by $T:M\longrightarrow\hat{M}$ where $M=\hat{M}\backslash H$ and $T=\hat{T}|_{M}$. The iterates of $T$ for $k=1,2,\dots$ are defined by $T^k:=(\hat{T}|_M)^k$ and we keep track of the iterates of points $x\in M$ as long as they do not enter the ``hole" $H$.\par
	For any point $x\in M$, the escape time is defined by the smallest natural number $n$ such that $\hat{T}^n(x)\in H$ and it is denoted by $\tau(x)$. We denote the set of all points that have not escaped after $n$ iterations by
	$$M^n=\{x\in M\ |\ \tau(x)>n\}.$$
	It is easy to see that
	$$M\supseteq M^1\supseteq M^2\supseteq\cdots.$$
	Let $\mu$ be a Borel probability measure on $M$. Then upper and lower bounds of escape rate of the measure $\mu$, respectively, $-\ln\overline{\lambda}$ and $-\ln\underline{\lambda}$ are defined by
	\begin{equation}\label{escaperate}
		\ln\overline{\lambda}:=\limsup_{n\rightarrow\infty}\frac{1}{n}\ln\mu(M^n)\ \ \ and\ \ \ \ln\underline{\lambda}:=\liminf_{n\rightarrow\infty}\frac{1}{n}\ln\mu(M^n).
	\end{equation}
	If $\ln\lambda:=\ln\underline{\lambda}=\ln\overline{\lambda}$, then the escape rate of $\mu$ equals $-\ln\lambda$.  
	Note that we assume $\ln0=-\infty$. If $\mu(M^n)$ is well defined, then it is called the survival probability after $n$ iterations. Clearly, the escape rate and survival probability depend on the measure $\mu$. For example, one can choose a Borel measure $\mu$ such that $\mu(M^n)=0$ for some $n$ or $\mu(M^n)=1$ for all $n$. If $\hat{T}$ is a measurable map and $H$ is a measurable set with respect to the Borel measure $\mu$, then the escape rate $\rho$ of $\mu$ is defined by,
	\begin{equation}\label{escaperate2}
		\rho:=-\lim_{n\rightarrow\infty}\frac{1}{n}\ln\mu(M^n),
	\end{equation}
	or equivalently,
	\begin{equation}\label{escaperate3}
		\rho=-\ln(\lim_{n\rightarrow\infty}\mu(M^n)^{\frac{1}{n}}).
	\end{equation}
	when $\lim_{n\rightarrow\infty}\mu(M^n)^{\frac{1}{n}}\neq0$, otherwise $\rho=\infty$.\par
	A Borel measure $\mu$ is called  conditionally invariant with respect to $T$ if
	\begin{equation}\label{condinv}
		\mu(A)=\frac{\mu(T^{-1}(A))}{\mu(T^{-1}(M))},
	\end{equation}
	for all Borel sets $A\subseteq M$. If $\mu$ is a conditionally invariant measure and $\lambda=\mu(M^1)$, then $-\ln\lambda$ is the escape rate of the open system with respect to the conditionally invariant measure $\mu$.\par
	The following constructive procedure ensures existence of an absolutely continuous conditionally invariant measure (a.c.c.i.m.) for an open dynamical system \cite{DY06} (an a.c.c.i.m. is a conditionally invariant measure that has density with respect to Lebesgue measure $m$). Assume that $\hat{T}:\hat{M}\longrightarrow\hat{M}$ admits a finite Markov partition $\mathcal{P}=\{E_1,\dots,E_k\}$ on $\hat{M}$.
	The existence of Markov partition assumes that $\hat{T}$ is a hyperbolic (i.e. a chaotic) dynamical system because Markov partitions are introduced and exist only for hyperbolic systems. Recall that the dynamical system is called hyperbolic if through almost every points of the phase space are passing stable and unstable manifolds.\par
	Let the hole $H$ be an element of Markov partition $\mathcal{P}$, and $T=\hat{T}|_M$ be the open dynamical system with hole $H$ (i.e. $M=\hat{M}\backslash H$). The (substochastic) transition matrix $P=[p_{ij}]$ of the partition $\mathcal{P}$ under the action of $T$ is defined by
	$$p_{ij}=\frac{m(E_i\cap T^{-1}(E_j))}{m(E_i)},$$
	where $1\leq i,j\leq k$ and $m$ is Lebesgue measure. Using refinements of the partition $\mathcal{P}$ given by $\mathcal{P}_n=\bigvee_{i=-n}^n\hat{T}^i(\mathcal{P})$, we obtain a sequence of substochastic matrices $P_n$. By Perron-Frobenius theorem, all $P_n$ have a positive leading (i.e. a maximal) eigenvalue $\lambda_n$. Under the usual conditions of aperiodicity and irreducibility on the matrices $P_n$, there exists an unique probability eigenvector $v_n$ corresponding to $\lambda_n$. If the Markov partition $\mathcal{P}$ is a generating partition, then the sequence of probability eigenvectors $v_n$ will converge to an a.c.c.i.m. $\mu$, where $\mu(M^1)=\lim_{n\rightarrow\infty}\lambda_n$. Hence, if we let $\lambda:=\mu(M^1)=\lim_{n\rightarrow\infty}\lambda_n$, then the escape rate of the open system $T:M\longrightarrow\hat{M}$ with respect to the a.c.c.i.m. $\mu$ is $-\ln\lambda$.
	%%%%%%%%%%%%%%%%%%%%%%%%%%%%%%%%%%%%%%%%%%%%%%%%%%%%%%%%%%%%%%%%%%%%%%%%%%%%%%%%%%%%%%%%%%%%%%%%%%%%%%%%%%%%%%%
	\section{Average of escape rates through elements of partition is larger than expected}
	Let $\hat{T}:\hat{M}\longrightarrow\hat{M}$ be a discrete-time dynamical system and $\mu$ be a ``natural" invariant probability measure on $\hat{M}$ with respect to $\hat{T}$ (i.e. $\mu(A)=\mu(\hat{T}^{-1}(A))$ for all Borel sets $A\subseteq\hat{M}$). We also assume that the ``closed" dynamical is ergodic with respect to $\mu$. The assumption of ergodicity ensures that almost all (i.e. measure one subset) orbits will eventually enter a ``hole" when the hole is a subset of positive measure of the phase space. We assume for simplicity that $\mu$ is absolutely continuous (i.e. has a density) with respect to Lebesgue measure.\par
	Consider a Markov partition $\{E_1,\dots,E_k\}$ of the phase space $\hat{M}$ with elements $E_i$ of positive measure (i.e. $\mu(E_i)>0$ for all $i$). Let $T_i=\hat{T}|_{\hat{M}\backslash E_i}$ be the open dynamical system with hole $E_i$. We assume $T_i$ admits an a.c.c.i.m. $\mu_i$ such that its escape rate $\rho_i$ is well-defined. We also denote the set of all points that have not escaped after $n$ iterations of $T_i$ by $M_i^n$. That means, 
	$$M_i^n:=\{x\in\hat{M}\backslash E_i\ |\ \tau_i(x)>n\}$$
	where $\tau_i(x)$ is the escape time of the open dynamical system $T_i$. Moreover, we assume that $T_i$ is a measurable map and $E_i$ is a measurable set with respect to the a.c.c.i.m. $\mu_i$. Under these assumptions for all $i=1,\dots,k$, we can set
	\begin{equation}\label{p_i}
		p_i:=\mu_i(M_i^1)=\lim_{n\rightarrow\infty}(\mu_i(M_i^n))^{1/n}.    
	\end{equation}
	Our goal is to prove the inequality
	\begin{equation}\label{ineq1}
		\langle\rho\rangle:=\sum_{i=1}^k\mu(E_i)\rho_i\geq-\ln(\sum_{i=1}^k\mu(E_i)p_i).
	\end{equation}
	If $p_i\neq0$ for all $i=1,\dots,k$, then by making use of (\ref{escaperate3}), we obtain
	$$\langle\rho\rangle=\sum_{i=1}^k\mu(E_i)\rho_i=-\sum_{i=1}^k\mu(E_i)\ln(p_i).$$
	Because $-\ln(x)$ is a convex function, the relation (\ref{ineq1}) follows from Jensen's inequality. Moreover, if $p_i=0$ for some $1\leq i\leq k$, then $\rho_i=\infty$, and it is easy to see that (\ref{ineq1}) is satisfied in this case as well.\par
	The left side of (\ref{ineq1}) is the average (with respect to the invariant probability measure $\mu$) of escape rates over all positions of the hole, i.e. over all elements of the Markov partition. The right hand side of (\ref{ineq1}) is equivalent to the escape rate of the system when we consider uniform density with a uniform leak from it where the hole's size is equal to the average of $1-\mu_i(M_i^1)$, $1\leq i\leq k$,  with respect to a probability measure $\mu$. Under these assumptions, the estimate of escape rate is 
	\begin{equation}\label{lb}
		-\ln(1-\sum_{i=1}^k\mu(E_i)(1-\mu_i(M_i^1))).
	\end{equation}
	Then, (\ref{p_i}) and (\ref{lb}) imply
	\begin{equation}\label{lb1}
		-\ln(1-\sum_{i=1}^k\mu(E_i)(1-\mu_i(M_i^1)))=-\ln(\sum_{i=1}^k\mu(E_i)\mu_i(M_i^1))=-\ln(\sum_{i=1}^k\mu(E_i)p_i).
	\end{equation}
	Hence, (\ref{ineq1}) and (\ref{lb1}) show that the (average) escape rate is faster than the estimate of escape rates when we assume that the density remains uniform in the system and the hole's size is the average of $1-\mu_i(M_i^1)$ with respect to the probability measure $\mu$. Remind that at each iteration of the map $T_i$ the same portion $1-\mu_i(M_i^1)$ of the remaining in the phase space points escapes through the hole $E_i$ with respect to the a.c.c.i.m. $\mu_i$. Therefore, the average of $1-\mu_i(M_i^1)$ is the average on proportional leaks of open systems $T_i$ with respect to their a.c.c.i.m. $\mu_i$.\par
	There are some other candidates for (naive) estimation of escape rate (see e.g. \cite{GDA12}). The first one is $-\ln(1-h)$ where $h=1/k$ is the average of holes' sizes (i.e. $h=\frac{\sum_{i=1}^k\mu(E_i)}{k}$) and $k$ is the number of elements in a Markov partition. We will denote this naive estimation with $N_1$. This naive comes from the assumption that the density of the a.c.c.i. measure remains uniform in the system and the hole’s size is equal to $h=1/k$.\par
	The second candidate $N_2$ for a naive estimate of escape rate was introduced in \cite{GDA12}. It is defined as $N_2:=-\sum_{i=1}^kh_i\ln(1-h_i)$ where $h_i=\mu(E_i)$. (Recall that $\mu$ is the invariant probability measure of the corresponding closed system). The naive estimate $N_2$ uses the assumption that the density of conditionally invariant measure is a constant in a complement of a corresponding hole. The estimate of escape rate of open system with hole $E_i$ is $-\ln(1-h_i)$ where $h_i=\mu(E_i)$. Then $N_2$ is the average of these estimates over different positions of the hole in the phase space (over the elements of Markov partition) with respect to probability measure $\mu$.\par
	By making use of the Jensen's inequality, we obtain
	\begin{equation}\label{N12}
		N_2=-\sum_{i=1}^kh_i\ln(1-h_i)\geq-\ln(\sum_{i=1}^kh_i(1-h_i))=-\ln(1-\sum_{i=1}^kh_i^2).
	\end{equation}
	It is easy to see that
	\begin{equation}\label{Var}
		\sum_{i=1}^kh_i^2\geq h=\frac{1}{k},
	\end{equation}
	since $1=\sum_{i=1}^k\mu(E_i)=\sum_{i=1}^kh_i$. Then, (\ref{N12}) and (\ref{Var}) imply
	$$N_2\geq-\ln(1-\sum_{i=1}^kh_i^2)\geq-\ln(1-h)=N_1.$$
	Therefore the naive estimation $N_2$ is greater than or equal to $N_1$. These two naive estimates become equal only  if elements of the partition have equal sizes (i.e. $\mu(E_i)=1/k$ for all $i$).
	%%%%%%%%%%%%%%%%%%%%%%%%%%%%%%%%%%%%%%%%%%%%%%%%%%%%%%%%%%%%%%%%%%%%%%%%%%%%%%%%%%%%%%%%%%%%%%%%%%%%%%%%%%%%%%%
	\section{Examples}
	%%%%%%%%%%%%%%%%%%%%%%%%%%%%%%%%%%%%%%%%%%%%%%%%%%%%%%%%%%%%%%
	\subsection{Skewed tent map}
	Consider a skewed tent map $T:[0,1]\longrightarrow[0,1]$ given by:
	\begin{equation}\label{tentmap}
		T(x)=\left\{
		\begin{aligned}
			& \frac{x}{x_0} \hspace{1cm} & x\in[0,x_0),\\
			& \frac{1-x}{1-x_0} & x\in[x_0,1],
		\end{aligned}
		\right.    
	\end{equation}
	where $x_0\in(0,1)$. Consider a Markov partition $\mathcal{E}=\{E_1,E_2\}$ of the phase space $[0,1]$, where $E_1=[0,x_0]$ and $E_2=[x_0,1]$. Let $\rho_i$ denote the escape rate of the open system with hole $E_i$. By making use of (\ref{escaperate2}) or (\ref{escaperate3}), we obtain
	$$\rho_1=-\ln(1-x_0), \hspace{2cm} \rho_2=-\ln(x_0).$$
	Thus,
	\begin{equation}\label{ex1-1}
		\langle\rho\rangle:=m(E_1)\rho_1+m(E_2)\rho_2=-x_0\ln(1-x_0)-(1-x_0)\ln(x_0),
	\end{equation}
	where $m$ is Lebesgue measure (the natural invariant probability measure of (\ref{tentmap})). From (\ref{p_i}), we get $p_1=1-x_0$ and $p_2=x_0$. Then,
	\begin{equation}\label{ex1-2}
		-\ln(\sum_{i=1}^2m(E_i)p_i)=-\ln(2x_0(1-x_0)).
	\end{equation}
	From (\ref{ex1-1}) and (\ref{ex1-2}) by using Jensen's inequality, we obtain
	\begin{equation}\label{ex1-3}
		\langle\rho\rangle\geq-\ln(\sum_{i=1}^2m(E_i)p_i).
	\end{equation}
	If $x_0=\frac{1}{2}$, then the equality will be satisfied in (\ref{ex1-3}).\par
	We compare now our estimate of the escape rate which is the lower bound (right hand side) of the inequality (\ref{ineq1}) with the naive estimate $N_1=-\ln(1-h)$ of escape rate, where $h$ is the average of holes' (elements of partition) sizes \cite{APT13,GDA12}. To do this we consider refinements of the partition $\mathcal{E}=\{[0,x_0],[x_0,1]\}$ under the skewed tent map $T$ to construct Markov partitions with $4,\ 8,\dots,\ 128,$ or $256$ elements.\par
	The Table \ref{table1} presents values of our estimation of escape rates of the skewed tent map (\ref{tentmap}), where $x_0=0.1,\ 0.2,\ 0.3,\ 0.4,$ or $0.5$ and the partition has $4,\ 8,\dots,\ 128,$ or $256$ elements. In Table \ref{table2}, naive estimations $N_1$ of the skewed tent map are presented, where the partition has $4,\ 8,\dots,\ 128,$ or $256$ elements. We know $N_1$ only depends on the number of elements in a partition, and it does not depend on the size of holes.\par
	\begin{table}[h]
		\begin{center}
			\begin{tabular}{ | c || c | c | c | c | c | c | c | } 
				\hline
				& 4 & 8 & 16 & 32 & 64 & 128 & 256\\
				\hline
				\hline
				0.1 & 0.77922 & 0.44239 & 0.28375 & 0.19638 & 0.14384 & 0.10949 & 0.08598\\
				\hline
				0.2 & 0.47400 & 0.25981 & 0.16685 & 0.11452 & 0.07286 & 0.04757 & 0.03175\\
				\hline
				0.3 & 0.37517 & 0.21720 & 0.11717 & 0.06234 & 0.03491 & 0.01987 & 0.01149\\
				\hline
				0.4 & 0.37047 & 0.17868 & 0.08023 & 0.03931 & 0.01990 & 0.01022 & 0.00529\\
				\hline
				0.5 & 0.42387 & 0.15808 & 0.06928 & 0.03297 & 0.01604 & 0.00792 & 0.00393\\
				\hline
			\end{tabular}
			\caption{The estimation of escape rate in (\ref{ineq1}) of the skewed tent map, where $x_0=0.1,\ 0.2,\ 0.3,\ 0.4,$ or $0.5$ and the partition has $4,\ 8,\dots,\ 128,$ or $256$ elements.}
			\label{table1}
		\end{center}
	\end{table}
	\begin{table}[h]
		\begin{center}
			\begin{tabular}{ | c || c | c | c | c | c | c | c | } 
				\hline
				& 4 & 8 & 16 & 32 & 64 & 128 & 256\\
				\hline
				\hline
				$N_1$ & 0.28768 & 0.13353 & 0.06453 & 0.03174 & 0.01574 & 0.00784 & 0.00391\\
				\hline
			\end{tabular}
			\caption{The naive estimation $N_1$ of escape rate of the skewed tent map, where the partition has $4,\ 8,\dots,\ 128,$ or $256$ elements.}
			\label{table2}
		\end{center}
	\end{table}
	Comparing the data of Table \ref{table1} with Table \ref{table2}, one can see our estimation of the escape rate is greater than the naive estimation $N_1$ for the skewed tent map.
	%%%%%%%%%%%%%%%%%%%%%%%%%%%%%%%%%%%%%%%%%%%%%%%%%%%%%%%%%%%%%%
	\subsection{Arnold's cat map}
	Let $\mathbb{T}^2$ be a torus, i.e the unit square $\mathbf{S}=[0,1]\times[0,1]$ with standard identification of its opposite sides. The Arnold's cat map $T:\mathbb{T}^2\longrightarrow\mathbb{T}^2$ is defined by
	$$T(\begin{bmatrix}x\\ y\end{bmatrix})=\begin{bmatrix}
	2 & 1\\
	1 & 1
	\end{bmatrix}\begin{bmatrix}x\\ y\end{bmatrix}\mod{1}.$$  
	The eigenvalues of $\begin{bmatrix}
	2 & 1\\
	1 & 1
	\end{bmatrix}$ are $\lambda_1=\frac{3+\sqrt{5}}{2}$ and $\lambda_2=\frac{3-\sqrt{5}}{2}$, where
	$$v_1=\begin{bmatrix}
	1\\
	\frac{\sqrt{5}-1}{2}
	\end{bmatrix},\ \ \ \ v_2=\begin{bmatrix}
	1\\
	-\frac{\sqrt{5}+1}{2}
	\end{bmatrix},$$
	are their corresponding eigenvectors, respectively.\par
\begin{figure}[h]
		\centering
		\includegraphics[width=10cm]{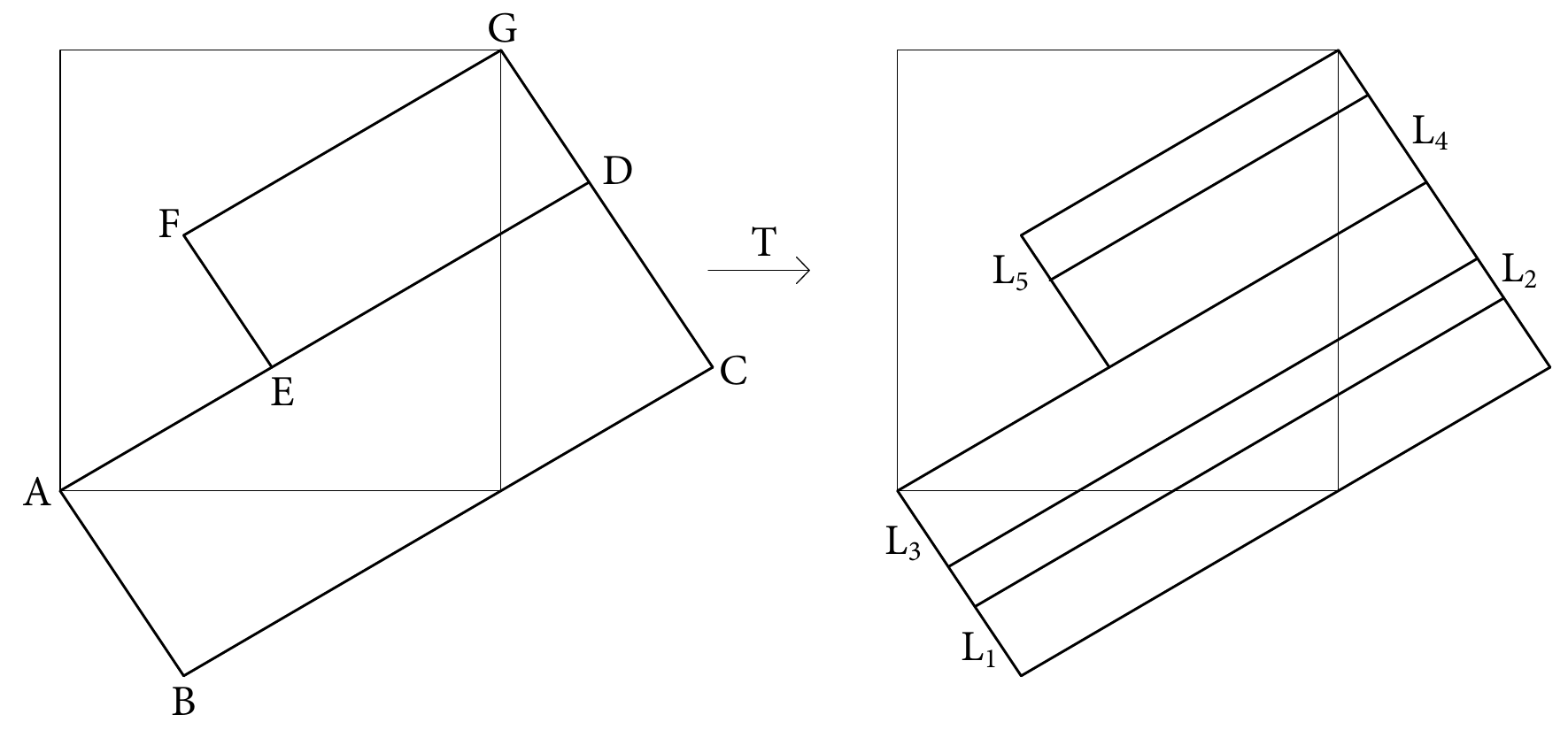}
		\caption{$B=(\frac{\sqrt{5}-1}{2\sqrt{5}},\frac{-1}{\sqrt{5}}), C=(\frac{\sqrt{5}+1}{\sqrt{5}},\frac{\sqrt{5}-1}{2\sqrt{5}}), D=(\frac{\sqrt{5}+3}{2\sqrt{5}},\frac{\sqrt{5}+1}{2\sqrt{5}}), E=(\frac{1}{\sqrt{5}},\frac{\sqrt{5}-1}{2\sqrt{5}}), F=(\frac{\sqrt{5}-1}{2\sqrt{5}},\frac{\sqrt{5}-1}{\sqrt{5}})$.}\label{fig1}
	\end{figure}

	It is easy to see that this map preserves the Lebesgue measure (area) on the torus (it is the natural invariant probability measure of Arnold's cat map). Consider rectangles $ABCD$ and $DEFG$ where their sides are parallel to the directions of eigenvectors $v_1$ and $v_2$ (as it is shown in Fig. \ref{fig1}). Take now the partition $\{R_1,R_2\}$ of $\mathbb{T}^2$, where $R_1$ ($R_2$) is the projection of rectangle $ABCD$ ($DEFG$) on the torus $\mathbb{T}^2$. Recall that the plane is a natural unfolding of the torus under identification of its parallel sides.\par
	The partition $\{R_1,R_2\}$ is not a generating one, but one can build the generating partition $\mathcal{L}=\{L_1,L_2,L_3,L_4,L_5\}$ from $\{R_1,R_2\}$ as follows \cite{AW70}:
	$$\begin{aligned}
	L_1\cup L_3 & =R_1\cap T(R_1)\\
	L_2 & =R_1\cap T(R_2)\\
	L_4 & =R_2\cap T(R_1)\\
	L_5 & =R_2\cap T(R_2)\\
	\end{aligned}$$
	See Fig. \ref{fig1} for more details.\par
	Let $T_0$ be the transition probability matrix of the partition $\mathcal{L}$ under the map $T$, then
	$$T_0=\begin{bmatrix}
	\frac{3-\sqrt{5}}{2} & 0 & \frac{3-\sqrt{5}}{2} & \sqrt{5}-2 & 0\\[6pt]
	\frac{3-\sqrt{5}}{2} & 0 & \frac{3-\sqrt{5}}{2} & \sqrt{5}-2 & 0\\[6pt]
	\frac{3-\sqrt{5}}{2} & 0 & \frac{3-\sqrt{5}}{2} & \sqrt{5}-2 & 0\\[6pt]
	0 & \frac{\sqrt{5}-1}{2} & 0 & 0 & \frac{3-\sqrt{5}}{2}\\[6pt]
	0 & \frac{\sqrt{5}-1}{2} & 0 & 0 & \frac{3-\sqrt{5}}{2}\\[6pt]
	\end{bmatrix}.$$
	Consider the open dynamical system built from Arnold's cat map with the hole $L_i$. If we denote the substochastic matrix of transition probabilities of this open system by $T_i$, then $T_i$ is equal to $T_0$ where the $i$th row of $T_0$ is replaced by zeros.\par
	Let $\mathcal{L}_n=\bigvee_{i=-n}^nT^i(\mathcal{L})$ denote the refinements of partition $\mathcal{L}$ under the action of map $T$ and $T_{i,n}$ be the substochastic transition matrix corresponding to the refined partition $\mathcal{L}_n$ of the open dynamical system with hole $L_i$. It is easy to check that $T_{i,n}$ and $T_i$ for all $n=1,2,\dots$ have the same leading eigenvalues $\lambda_i$. Also, we know $\mathcal{L}$ is a generating Markov partition and all $T_{i,n}$ satisfy  aperiodicity and irreducibility conditions. Therefore, there is an a.c.c.i.m. $\mu_i$ where the escape rate of open system with hole $L_i$ with respect to $\mu_i$ is $-\ln\lambda_i$. Here, we have
	$$\lambda_1=\lambda_2=\lambda_3=\lambda_4=3-\sqrt{5},\ \ \ \lambda_5=\frac{1+\sqrt{2}}{2}(3-\sqrt{5}).$$
	Hence, the escape rates corresponding to the a.c.c.i. measures $\mu_i$ are 
	$$\rho_1=\rho_2=\rho_3=\rho_4=-\ln(3-\sqrt{5}),\ \ \ \rho_5=-\ln(\frac{1+\sqrt{2}}{2}(3-\sqrt{5})).$$
	Thus the average of escape rates with respect to Lebesgue measure is
	\begin{equation}\label{ex2-1}
		\langle\rho\rangle=\sum_{i=1}^5m(L_i)\rho_i\simeq0.2494,
	\end{equation}
	where
	$$m(L_1)=m(L_3)=\frac{3-\sqrt{5}}{2}\sqrt{\frac{\sqrt{5}+3}{10}},\ \ \ m(L_2)=(\sqrt{5}-2)\sqrt{\frac{\sqrt{5}+3}{10}},$$
	$$m(L_4)=\frac{\sqrt{5}-1}{2}\left(1-\sqrt{\frac{\sqrt{5}+3}{10}}\right),\ \ \ m(L_5)=\frac{3-\sqrt{5}}{2}\left(1-\sqrt{\frac{\sqrt{5}+3}{10}}\right).$$
	From (\ref{p_i}), we obtain
	$$p_1=p_2=p_3=p_4=\lambda_1=3-\sqrt{5},\hspace{1.5cm} p_5=\lambda_5=\frac{1+\sqrt{2}}{2}(3-\sqrt{5}).$$
	Therefore,
	\begin{equation}\label{ex2-2}
		-\ln(\sum_{i=1}^5m(L_i)p_i)\simeq0.2476.
	\end{equation}
	By comparing (\ref{ex2-1}) and (\ref{ex2-2}), we see that the average of escape rates over the elements of Markov partition $\mathcal{L}$ is greater than the estimation of escape rate for the Arnold's cat map.
	%%%%%%%%%%%%%%%%%%%%%%%%%%%%%%%%%%%%%%%%%%%%%%%%%%%%%%%%%%%%%%
	\subsection{The Ulam$-$von Neumann logistic map}
	There are only a few examples of nonlinear systems where the escape rate is studied. Consider the nonlinear Ulam$-$von Neumann map of the unit interval \cite{ulam},
	\begin{equation}\label{logistic}
		T(x)=4x(1-x),\hspace{1.5cm}x\in[0,1].
	\end{equation}
	It is well-known that this map is metrically conjugate to the tent map (defined by (\ref{tentmap}) when $x_0=\frac{1}{2}$), where the conjugate map $U$ is given by,
	$$U(x)=\sin^2(\frac{\pi x}{2}).$$
	The invariant probability measure $\mu$ of Ulam$-$von Neumann logistic map (\ref{logistic}) has non-uniform density $f(x)$ with respect to Lebesgue measure $m$, where
	$$f(x)=\frac{1}{\pi\sqrt{x(1-x)}}.$$
	We consider the following Markov partition on the unit interval of the logistic map,
	$$\mathcal{P}=\left\{P_i\ |\ P_i=\left[\sin^2\left(\frac{i\pi}{2^{n+1}}\right),\sin^2\left(\frac{(i+1)\pi}{2^{n+1}}\right)\right],\ i=0,1,2,\dots,2^n-1\right\},$$
	where the partition $\mathcal{P}$ is the image of the natural Markov partition 
	$$\mathcal{E}=\left\{E_i\ |\ E_i=\left[\frac{i}{2^n},\frac{i+1}{2^n}\right],\ i=0,1,2,\dots,2^n-1\right\},$$
	of tent map under the conjugate map $U(x)$. It is easy to see that $\mu(P_i)=m(E_i)=\frac{1}{2^n}$ for all $i=0,1,\dots,2^n-1$. Moreover, the transition probability matrix of the partition $\mathcal{P}$ under Ulam$-$von Neumann logistic map is the same as the transition matrix of the tent map for the partition $\mathcal{E}$. These two maps also have the same transition matrices on the corresponding refinements of their partitions.\par
	Let $T_i$ denote the open dynamical system of the logistic map with hole $P_i$. Clearly, $T_i$ will have the same substochastic transition matrices as the open tent map with hole $E_i$ on refinements of their corresponding Markov partitions. Hence, the escape rate $\rho_{T_i}$ of $T_i$ is equal to the escape rate $\rho_i$ of open tent map with hole $E_i$. It implies
	$$\sum_{i=0}^{2^n-1}\mu(P_i)\rho_{T_i}=\sum_{i=0}^{2^n-1}m(E_i)\rho_i.$$
	For the same reason as before, the lower bound (the right hand side) of the inequality (\ref{ineq1}) (i.e. our naive estimate of the escape rate) will be the same for both these systems. Thus, the relation between our estimate of the escape rate (lower bound of (\ref{ineq1})) and the naive estimate $N_1$ will be valid in this case as well (see Tables \ref{table1} and \ref{table2}).
	%%%%%%%%%%%%%%%%%%%%%%%%%%%%%%%%%%%%%%%%%%%%%%%%%%%%%%%%%%%%%%%%%%%%%%%%%%%%%%%%%%%%%%%%%%%%%%%%%%%%%%%%%%%%%%%
	\section{Concluding remarks}
	We have shown that for chaotic maps which admit a finite generating Markov partition, the averaged over the elements of the MP escape rate exceed a naive estimate of the escape rate.\par
	A natural question would be to analyze relations between these and other possible estimates of a global (average) escape rate in nonlinear dynamical systems. We believe that our argument based on convexity will be still an important basic tool in this case as well.\par
	A standard approach going back to Sinai is to consider a sequence of Markov partitions with smaller and smaller elements \cite{Sin72}. Then the map becomes closer and closer to a linear on the elements of Markov partitions, and thus the entire dynamical system is approximated by a sequence of Markov chains. To perform such proofs, it seems that the higher (second) order approximation for escape rate in terms of the size of a ``hole" obtained in \cite{GDA12} could be quite useful. It is also worthwhile to mention in this respect that it was observed numerically \cite{PP} that in (nonlinear) logistic maps the process of escape also slows down near periodic orbits.\par
	There is also a direct connection between Markov chains, dynamical systems with Markov partitions, and transport (dynamics) on networks \cite{AB,BB}. Namely, the adjacency matrix of a network (which has the entry $1$ if the element $i$ is connected by an edge to the element $j$, or otherwise $0$ entries), can always be considered as a structural matrix of a Markov chain or as a transition matrix of a topological Markov chain (directed graph). Therefore, the results about open Markov systems are applicable to networks with leaking elements and allow for estimates of average leaking, most leaking sites, etc \cite{AB10,BB}.
	%%%%%%%%%%%%%%%%%%%%%%%%%%%%%%%%%%%%%%%%%%%%%%%%%%%%%%%%%%%%%%%%%%%%%%%%%%%%%%%%%%%%%%%%%%%%%%%%%%%%%%%%%%%%%%%
	
\end{document}